\documentclass[10pt,leqno,amssymbols]{article}

\usepackage{amsmath, amsthm, amsfonts, amssymb}
\usepackage[english]{babel}
\usepackage{amsmath}
\usepackage{amsfonts}
\usepackage{amssymb}
\usepackage{graphicx,color}

\newtheorem{thm}{Theorem}[section]
\newtheorem{cor}[thm]{Corollary}
\newtheorem{lem}[thm]{Lemma}
\newtheorem{prop}[thm]{Proposition}
\theoremstyle{definition}
\newtheorem{defn}[thm]{Definition}
\theoremstyle{remark}
\newtheorem{rem}[thm]{Remark}
\numberwithin{equation}{section}

\title{\textbf{ Norden structures on cotangent bundles}}

\author{\Large {Antonella Nannicini}}

\date{ }

\begin{document}
\maketitle

\noindent{ \large {\bf{Abstract.}}} 
We study prolongation of Norden structures on manifolds to their generalized tangent bundles and to their cotangent bundles. In particular, by using methods of generalized geometry,  we prove that the cotangent bundle of a complex Norden manifold $(M,J,g)$  admits a structure of  Norden manifold, $(T^{\star}(M),\tilde J, \tilde g)$. Moreover if $(M,J,g)$ has flat natural canonical connection then $\tilde J$ is integrable, that is $(T^{\star}(M),\tilde J, \tilde g)$ is a complex Norden manifold. Finally we prove that  if $(M,J,g)$ is  K\"ahler Norden flat then  $(T^{\star}(M),\tilde J, \tilde g)$ is  K\"ahler Norden flat.\\


\section{Introduction}

Let $M$  be a smooth manifold, let $ T(M)$ be the tangent bundle and let $ T^{\star}(M)$ be the cotangent bundle of $M$.  $ E=T(M) \oplus T^{\star}(M)$
 is called  the  \emph{generalized tangent bundle} of $M$. In \cite{H} Hitchin introduced the concept of generalized complex structure, he consided  complex structures of the generalized tangent bundle compatible with the standard metric of neutral signature on $E$. In previous papers, \cite{N1}, \cite{N2}, \cite{N3}, \cite{N4}, \cite{N5}, \cite{N6}, we defined and studied a class of complex structures on the generalized tangent bundle compatible with the standard symplectic structure of $E$, that is  a class of  \emph{pseudo calibrated generalized complex structures}. In particular in \cite{N5} we described generalized complex structures defined naturally by Norden structures. In this paper we study prolongation of Norden structures on manifolds to their generalized tangent bundles and to their cotangent bundles. Precisely, by using methods of generalized geometry,  we prove that the cotangent bundle of a complex Norden manifold, $(M,J,g)$,  admits a structure of  Norden manifold, $(T^{\star}(M), \tilde J, \tilde g)$. Moreover if $(M,J,g)$ has flat natural canonical connection then $\tilde J$ is integrable, that is $(T^{\star}(M),\tilde J, \tilde g)$ is a complex Norden manifold. Finally we prove that if $(M,J,g)$ is  K\"ahler Norden flat then  $(T^{\star}(M),\tilde J, \tilde g)$ is  K\"ahler Norden flat.\\
The paper is organized as follows. In sections 2 we introduce basics on generalized tangent bundles, generalized complex structures and Norden manifolds. Original results are contained in section 3, were  we construct  Norden structures on generalized tangent bundles, and in section 4,  where we construct  Norden structures on cotangent bundles and we describe curvature properties. \\

\noindent {\bf {Acknowledgements.}} The author's research was partially supported by the following grants of the Italian Ministry of Education (MIUR): PRIN Variet\`a reali e complesse (2010NNBZ78) and by GNSAGA of INDAM.

\section{Preliminaries}
\subsection{Generalized tangent bundle}
\noindent  Let $M$ be a smooth manifold of real dimension $n$ and let $ E=T(M) \oplus T^{\star}(M)$
 be the generalized tangent bundle of $M$. Smooth sections of  $E$ are elements $X+\xi \in C^\infty (E) $ where $X \in C^\infty(T(M))$ is a vector field and $\xi \in C^\infty( T^{\star}(M))$ is a $1$-form. $E$ is equipped with a natural symplectic structure, $(\, , \,)$, defined on two elements  $X+\xi, \, Y+\eta \in C^\infty (E)$ by:
\begin{equation}
(X+\xi,Y+\eta)=-{1\over 2}(\xi (Y)- \eta (X)).\label{10}
\end{equation}
\noindent Moreover $E$ is equipped with a natural indefinite metric, $<\, , \,>$, defined by:
\begin{equation}
<X+\xi,Y+\eta>=-{1\over 2}(\xi (Y)+ \eta (X)); \label{20}
\end{equation}
$<\, ,\,>$ is non degenerate and of signature $(n,n)$. \\
A linear connection on $M$, $\nabla$, defines a bracket on $C^ \infty (E)$, $[\, ,\,]_\nabla$, as follows:
\begin{equation}
[X+\xi,Y+\eta]_{\nabla} = [X,Y]+ {\nabla}_{X} {\eta} - {\nabla}_{Y} {\xi} \label{30}
\end{equation}
\noindent where $[\, ,\,]$ means Lie bracket of vector fields. \\

\begin{lem} (\cite{N1}) For all $X,Y \in C^ \infty (T(M))$, for all $\xi , \eta \in C^ \infty (T^{\star}(M))$ and for all $f\in C^ \infty (M)$ we have:\\
\noindent 1. $[X+\xi,Y+\eta]_{\nabla} =-[Y+\eta,X+\xi]_{\nabla} $ \\
\noindent 2. $[f(X+\xi),Y+\eta]_{\nabla} =f[X+\xi,Y+\eta]_{\nabla} -Y(f)(X+\xi)$  \\
\noindent 3. Jacobi's identity holds for $[\, ,\,]_{\nabla}$ if and only if $\nabla$ has zero curvature.
\end{lem}

\subsection{Pseudo calibrated generalized complex structures}

\begin{defn}
A \emph {generalized complex structure} on $M$ is an endomorphism $\widehat {J} : E \rightarrow E$ such that $ {\widehat {J}}^2=-I.$
\end{defn}

\begin{defn}
A generalized complex structure $\widehat{J}$ is called \emph {pseudo calibrated} if it is  $ (\, ,\, )$-invariant and if the bilinear symmetric form defined by $(\, ,\widehat {J})$ on $T(M)$ is non degenerate. Moreover  $\widehat{J}$ is called \emph {calibrated} if it is pseudo calibrated and  $(\, ,\widehat {J})$ is positive definite.
\end{defn}

\noindent From the definition we obtain that a generalized pseudo calibrated complex structure has the following block matrix form: 
\begin{equation}
\widehat {J} =  \begin{pmatrix} J & -(I+J^2)g^{-1} \cr 
g & -J^{\star}\cr
\end{pmatrix}\label{40}
\end{equation} 

\noindent where $g:T(M) \rightarrow T^{\star}(M)$ is identified to the bemolle musical isomorphism of a pseudo Riemannian metric $g$ on $M$, $J:T(M) \rightarrow T(M)$ is a $g$-symmetric operator and $J^{\star}:T^{\star}(M) \rightarrow T^{\star}(M)$ is the dual operator of $J$ defined by $J^{\star}(\xi)(X)= \xi (J(X))$.\\
\noindent $\widehat {J}$ is calibrated if and only if  $g$ is a Riemannian metric, namely:
\begin{equation}
(g(X))(Y)=g(X,Y)=2(X, \widehat {J} Y).\label{45}
\end{equation}

\subsection{Integrability}

\begin{lem} (\cite{N2}) Let $\widehat {J} : E \rightarrow E$ be a generalized complex structure on $M$ and let 

\begin{equation} {N^{\nabla}({\widehat {J}} ):C^{\infty}(E) \times C^{\infty}(E) \rightarrow C^{\infty}(E)}\label{50}
\end{equation}
defined for all $\sigma, \tau \in C^{\infty}(E)$  by:

\begin{equation} N^{\nabla}(\widehat {J})(\sigma,\tau)=[\widehat {J} \sigma,\widehat {J} \tau]_{\nabla}-\widehat {J}[\widehat {J} \sigma, \tau]_{\nabla} - \widehat {J}[\sigma,\widehat {J} \tau]_{\nabla} - [\sigma, \tau]_{\nabla}.\label{60}
\end{equation} 
$N^{\nabla}(\widehat {J})$  is a skew symmetric tensor called the {\bf {Nijenhuis tensor}} of  $\widehat {J}$ with respect to $\nabla$.
\end{lem}

\noindent  Let $E^{\mathbb{C}}=(T(M) \oplus T^{\star}(M)) \otimes {\mathbb{C}}$ be the complexified generalized tangent bundle. The splitting in $\pm i$ eigenspaces of $\widehat J$ is denoted by: $E^{\mathbb{C}}=E^{1,0}_{\widehat J} \oplus E^{0,1}_{\widehat J}$, with $E^{0,1}_{\widehat J}={\overline {{E^{1,0}_{\widehat J}}}}.$  Let $P_+:E^{\mathbb{C}} \rightarrow E^{1,0}_{\widehat J}$ and $P_-:E^{\mathbb{C}} \rightarrow E^{0,1}_{\widehat J}$ be the projection operators: $P_{\pm}={1\over {2}} (I\mp i{\widehat J})$. The following holds.\\

\begin{lem} (\cite{N5}) For all $\sigma, \tau \in C^{\infty}(E^{\mathbb{C}})$  we have:
\begin{equation} P_{\mp} [P_{\pm}(\sigma),P_{\pm}(\tau)]_{\nabla}=-{1\over 4}P_{\mp}(N^{\nabla}(\widehat {J})(\sigma,\tau)).\label{70}
\end{equation}
\end{lem}

\begin{cor} For any linear connection $\nabla$ on $M$ $E^{1,0}_{\widehat J}$ and $ E^{0,1}_{\widehat J}$ are $[\, , \,]_{\nabla}$-involutive if and only if $N^{\nabla}(\widehat {J})=0$.
\end{cor}

\begin{defn}
Let  $\widehat {J} : E \rightarrow E$ be a generalized complex structure on $M$, $\widehat {J}$ is called  \emph {$\nabla $-integrable} if   $N^{\nabla}(\widehat {J})=0$.
\end{defn}

\noindent Let $(M,g)$ be a pseudo Riemannian manifold, let $\nabla$ be a linear connection on $M$, the\emph{ torsion} of $\nabla$, $T^{\nabla}$, and the\emph { exterior differential} associated to $\nabla$ acting on $g$, $(d^{\nabla}g)$, are defined on $X,Y \in C^{\infty}(T(M)))$ respectively by:
\begin{equation} T^{\nabla}(X,Y)={\nabla}_X Y-{\nabla}_Y X -[X,Y] \label{74}
\end{equation}
\begin{equation} (d^{\nabla}g)(X,Y)=({\nabla}_X g)(Y)-({\nabla}_Y g)(X)+g(T^{\nabla}(X,Y)).\label{77}
\end{equation}
\noindent Let $J$ be a $g$-symmetric operator on $T(M)$ and let $N(J)$ be the Nijenhuis tensor of $J$, defined on  $X,Y \in C^{\infty}(T(M))$ by:
\begin{equation} {N(J)(X,Y)=[JX,JY]-H[JX,Y]-J[X,JY]+J^2[X,Y].}\label{80}
\end{equation}
Let us suppose $J^2=- I$.  Let $\widehat{J}$ be the pseudo calibrated generalized complex structure defined by $g$ and $J$: 
\begin{equation} \widehat {J} =  \begin{pmatrix} J & 0 \cr 
g & -J^{\star} \cr
\end{pmatrix}. \label{90}
\end{equation}
\noindent $\nabla$ - integrability of $\widehat{J}$ is described as follows.\\

\begin{thm} (\cite{N5}, \cite{N6}) The pseudo calibrated generalized complex structure\\
$\widehat {J}= \begin{pmatrix} J & 0 \cr 
g & -J^{\star} \cr
\end{pmatrix}$ is $\nabla$ -integrable if and only  if for all $X,Y \in C^{\infty}(T(M))$ the following conditions hold: 
\begin{equation} \left \{ \begin{array}{l} N(J) = 0\\

 {\nabla}J=0\\ 

(d^{\nabla}g)(JX,Y)+(d^{\nabla}g)(X,JY)=0.
\end{array} 
\right.\label{100}
\end{equation}
\end{thm}

\subsection{Norden manifolds}

Norden manifolds were introduced by  Norden in \cite{N} and then studied also as almost complex manifolds with B-metric and anti K\"ahlerian manifolds, \cite{BFV}, \cite{GM}. They have applications in mathematics and in theoretical physics.  \\

\begin{defn}
Let $(M,J)$ be an almost complex manifold and let $g$ be a pseudo Riemannian metric on $M$ such that $J$ is a $g$-symmetric operator, $g$ is called \emph {Norden metric} and $(M,J,g)$ is called \emph {Norden manifold}.
\end{defn}

\begin{defn}
Let $(M,J,g)$ be a Norden manifold, if $J$ is integrable then $(M,J,g)$ is called \emph {complex Norden manifold}.  
\end{defn}

\noindent The following result is well known.

\begin{thm} (\cite{GM})
Let $(M,J,g)$ be a complex Norden manifold, there exists a unique linear connection $D$ with torsion $T$ on $M$ such that:
\begin{equation}
\left \{ \begin{array}{l} (D_Xg)(Y,Z)=0\\ 

T(JX,Y)+T(X,JY)=0\\ 

g(T(X,Y),Z)+g(T(Y,Z),X)+g(T(Z,X),Y)=0.
\end{array}
\right.\label{110}
\end{equation}
for all $X,Y,Z \in C^{\infty}(T(M))$. $D$ is called the \bf {natural canonical connection}.
\end{thm}

\begin{rem}
The natural canonical connection $D$ is defined by: 
\begin{equation} D_X Y={\nabla}_X Y- {1\over 2}J({\nabla}_X J)Y\label{113}
\end{equation}
where $\nabla $ is the Levi Civita connection of $g$. In particular $DJ=0$.
\end{rem}

\noindent Then, from Theorem 2.8 and Theorem 2.11,  we get the following.\\

\begin{prop} (\cite{N5})
Let $(M,J,g)$  be a complex Norden manifold and let $D$ be the natural canonical connection on $M$, the generalized complex structure on $M$  defined by $J$ and $g$:
\begin{equation}{\widehat {J}= \begin{pmatrix}J & 0\cr 
g &-J^{\star} \cr
\end{pmatrix}}\label{120}
\end{equation}
\emph   {is} $D$ -integrable.
\end{prop}

\begin{defn} Let $(M,J,g)$ be a Norden manifold and let $\nabla$ be the Levi Civita connection of $g$, $(M,J,g)$ is called \emph {K\"ahler Norden manifold} if $ \nabla J=0$. 
\end{defn}

\begin{rem}
For a K\"ahler Norden manifold  $(M,J,g)$ the structure $J$ is integrable and the natural canonical connection is the Levi Civita connection of $g$.
\end{rem}

\section{Norden structures on generalized tangent bundles}

\subsection{Norden metric on $T(M) \oplus T^{\star}(M)$}

\noindent Let $(M,J,g)$ be a Norden manifold and let $ E=T(M) \oplus T^{\star}(M)$  be the generalized tangent bundle of $M$. Let $X+\xi, \, Y+\eta$ be smooth sections of $E$ and let $\natural : T^{\star}(M) \rightarrow T(M)$ be the inverse of the bemolle musical isomorphism of $g$, $\flat : T(M) \rightarrow T^{\star}(M)$, $\flat (X)(Y)=g(X,Y)$, for all $X, \, Y \in C^{\infty}(T(M))$.\\

\noindent  We define:
\begin{equation}\hat g(X+\xi, Y+\eta)=g(X,Y)+{1\over 2}g(JX,\natural \eta)+{1\over 2}g(\natural \xi, JY)+g(\natural \xi, \natural \eta).\label{130}
\end{equation}

\noindent We have the following.

\begin{prop} $\hat g$ is symmetric and the pseudo calibrated generalized complex structure on $M$ defined by $g$ and $J$,  $\widehat {J} =  \begin{pmatrix} J & 0 \cr 
g & -J^{\star} \cr
\end{pmatrix}$, is a $\hat g$-symmetric operator.
\end{prop}

\noindent {\bf{Proof.}} The symmetry of $\hat g$ follows immediately from the definition. In order to verify that $\widehat {J} $  is a $\hat g$-symmetric operator let us compute: \\

$\hat g(\hat J (X+\xi), Y+\eta)=g(JX+g(X)-J^{\star}(\xi),Y+\eta)$\\

$=g(JX,Y)+{1\over 2}g(-X,\natural \eta)+{1\over 2}g(X-\natural J^{\star}(\xi), JY)+g(X-\natural J^{\star}(\xi), \natural \eta)$\\

$=g(JX,Y)+{1\over 2}g(X,\natural \eta)+{1\over 2}g(X,JY)-{1\over 2}g(\natural J^{\star}(\xi),JY)-g(\natural J^{\star}(\xi), \natural \eta).$\\

\noindent On the other hand we have:\\

$\hat g(X+\xi,\hat J (Y+\eta))=g(X+\xi,JY+g(Y)-J^{\star}(\eta))$\\

$=g(X,JY)+{1\over 2}g(JX,Y-\natural J^{\star}( \eta))+{1\over 2}g(\natural \xi, -Y)+g(\natural \xi,Y- \natural J^{\star}(\eta))$\\

$=g(X,JY)+{1\over 2}g(JX,Y)-{1\over 2}g(JX,\natural J^{\star}(\eta))+{1\over 2}g(\natural \xi,Y)-g(\natural \xi,\natural J^{\star}(\eta)).$\\

\noindent From the property $J^{\star}=\flat(J \natural)$  and from the fact that $(M,J,g)$ is a Norden manifold we get the statement. $\square$\\

\begin{defn} Let $(M,J,g)$ be a Norden manifold, $(\hat J, \hat g)$ defined as before is called \emph {Norden structure on  $T(M) \oplus T^{\star}(M)$} \emph{defined by} $(J,g)$.
\end{defn}

\subsection{Natural canonical connection on $T(M) \oplus T^{\star}(M)$}

 {Let} $(M,J,g)$  be a complex Norden manifold and let $D$ be the natural canonical connection on $M$, we define

\begin{equation} \hat D : C^{\infty}(T(M)) \times C^{\infty}(T(M)\oplus T^\star (M))\rightarrow C^{\infty}(T(M) \oplus T^{\star}(M)) \label{133}
\end{equation}
by:

\begin{equation}{\hat D}_X (Y+\eta)={D_X} Y+{D_X} {\eta}  \label{132}
\end{equation}

\noindent  for all $ X\in  C^{\infty}(T(M)), \,  Y+\eta\in C^{\infty}(T(M)\oplus T^\star (M))$.\\

\noindent We have:\\
\begin{prop} The following condition hold:
\begin{equation}\left \{ \begin{array}{l}{\hat D}_X{\hat J}=0\\ 

\vspace{0.1cm}
{\hat D}_X{\hat g}=0.

\end{array}
\right.\label{135}
\end{equation}

\noindent for all $X \in C^{\infty}(T(M))$.
\noindent Moreover if $D$ is flat then $\hat D$ is flat.\\
$\hat D$ is called the {\bf natural canonical connection} of $T(M)\oplus T^\star (M)$.\\

\end{prop}

\noindent {\bf{Proof.}} From the definition of $\hat J$ and from the properties of $D$ we get:\\

${\hat D}_X {\hat J}(Y+ \eta)={\hat D}_X (JY+g(X)-J^{\star}\eta)=D_X JY+D_X(g(Y)-J^{\star} \eta)$\\

$=JD_X Y+g(D_X Y)-J^{\star}(D_X \eta)=\hat J(D_X Y)+\hat J(D_X \eta)=\hat J({\hat D}_X (Y+\eta))$\\

\noindent for all  $ X\in  C^{\infty}(M), \,  Y+\eta\in C^{\infty}(T(M)\oplus T^\star (M))$.\\

\noindent Moreover from the definition of $\hat g$ we get:\\

$X\hat g(Y+\eta,Z+\zeta)-\hat g({\hat D}_X(Y+\eta),Z+\zeta)-\hat g(Y+\eta, {\hat D}_X (Z+\zeta))$\\

$=X\{g(Y,Z)+{1\over 2}g(JY,\natural \zeta)+{1\over 2}g(\natural \eta,JZ)+g(\natural \eta,\natural \zeta)\}+$\\

$-\{g({D}_X Y,Z)+{1\over 2}g(J(D_X Y),\natural \zeta)+{1\over 2}g(\natural(D_X \eta),JZ)+g(\natural(D_X\eta),\natural \zeta)\}+$\\

$-\{g(Y,{D}_XZ)+{1\over 2}g(J Y,\natural (D_X \zeta)+{1\over 2}g(\natural \eta,J(D_X Z))+g(\natural \eta,\natural (D_X \zeta))\}$\\

$=0$\\

\noindent for all $X, Y, Z \in C^{\infty}(T(M)),  \, Y+\eta, \,  Z+\zeta \in C^{\infty}(T(M)\oplus T^\star (M)).\\$

\noindent Finally,  if $R^{D}$ denotes the curvature tensor of $D$:\\

$R^D (X,Y)=D_X D_Y-D_Y D_X- D_{[X,Y]} $\\

\noindent and $R^{\hat D}$ the curvature tensor of $\hat D$:\\

$R^{\hat D}(X,Y)(Z+\zeta)={\hat D}_X{\hat D}_Y (Z+\zeta)-{\hat D}_Y{\hat D}_X(Z+\zeta)-{\hat D}_{[X,Y]}(Z+\zeta),$\\

\noindent  we have:\\

$R^{\hat D}(X,Y)(Z+\zeta)=R^D(X,Y)Z+R^D(X,Y)(\zeta)$\\

\noindent for all $X, Y,Z  \in C^{\infty}(T(M)),  \, \zeta \in C^{\infty}( T^\star (M)).$ \\

\noindent Then the proof is complete. $\square$\\

\section{Norden structures on cotangent bundles}

\subsection{Complex Norden structure on $T^{\star}(M)$}
\noindent Let $(M,J,g)$ be a Norden manifold and let $\nabla$ be a linear connection on $M$. $\nabla$ defines the decomposition in horizontal and vertical subbundles of $T(T^\star(M))$:
\begin{equation} T(T^{\star}(M))=T^H(T^{\star}(M))\oplus T^V(T^{\star}(M)).\label{140}
\end{equation}
Let $\left\{ x^{1},...,x^{n}\right\} $ be local coordinates on $M$, let $\left\{ {\tilde{x}}^{1},..., {\tilde{x}}^{n},y_1,...,y_n\right\} $ be the corresponding local coordinates on $T^{\star}(M)$ and let  $\left\{ \dfrac{\partial }{\partial {{\tilde{x}}^1}},...,\dfrac{\partial }{\partial  {\tilde{x}}^{n}},\dfrac{\partial }{\partial {y_{1}}},...,\dfrac{\partial }{\partial {y_{n}}}\right \} $ be local frames
on $T(T^{\star}(M))$. We have:
\begin{equation}
\left( \dfrac{\partial }{\partial
{\tilde x}^{i}}\right)^H={\dfrac{\partial }{\partial
{\tilde x}^{i}}}+y_k{\Gamma}^k_{il}{\dfrac{\partial }{\partial
y_{l}}} \label{142}
\end{equation}
\begin{equation}\left( \dfrac{\partial }{\partial
{\tilde x}^{i}}\right)^V=-y_k{\Gamma}^k_{il}{\dfrac{\partial }{\partial
y_{l}}}\label{144}
\end{equation}
\begin{equation}
\left( \dfrac{\partial }{\partial
{y_{i}}}\right)^H=0\label{146}
\end{equation}
\begin{equation}
\left ( \dfrac{\partial }{\partial
{y_{i}}}\right)^V= \dfrac{\partial }{\partial
{y_{i}}}
\label{150}
\end{equation}
\noindent where $i,k,l$ run from $1$ to $n$, $\Gamma^k_{il}$ are Christoffel's symbols of $\nabla$ and we used Einstein's convention on repeated indices.\\
\noindent Let us denote  $X_i= \dfrac{\partial }{\partial
{\tilde x}^{i}}$, we have:\\
\begin{equation}
\left[X_i^H,X_j^H\right]=y_k{R}^k_{ijl}{\dfrac{\partial }{\partial
y_{l}}}\label{158}
\end{equation}
\begin{equation}
\left[X_i^H,{\dfrac{\partial }{\partial
y_{j}}}\right]=-{\Gamma}^j_{il}{\dfrac{\partial }{\partial
y_{l}}},
\label{160}
\end{equation}
\noindent where $R$ is the curvature tensor of $\nabla$: 
\begin{equation}
R(X_i,X_j)X_l=R^k_{ijl}X_k=\left(\nabla_{X_i}\nabla_{X_j}-\nabla_{X_j}\nabla_{X_i}-\nabla_{[X_i,X_j]}\right)X_l.
\end{equation}

\noindent We recall the following.\\

\begin{prop} (\cite{N2})
Let $M$ be a smooth manifold and let $\nabla$  be a linear connection on $M$, there is a bundle morphism:
\begin{equation}
{\Phi}^{\nabla} :T(M) \oplus T^{\star}(M) \rightarrow T(T^{\star}(M))\label{170}
\end{equation} 
which is an isomorphism on the fibres and such that \\
\noindent 1. $({\Phi}^{\nabla})^{\star}(\Omega)=-2(\, , \,) $ if and only if $\nabla$  has zero torsion, \\
\noindent 2. $({\Phi}^{\nabla})( [\, ,\,]_{\nabla})=[{\Phi}^{\nabla},{\Phi}^{\nabla}] $  if and only if $\nabla$ has zero curvature, \\ 
\noindent where $\Omega$  is the canonical symplectic form on $T^{\star}(M)$ defined by the Liouville $1$-form.
\end{prop}

\noindent In local coordinates we have the following expressions:

\begin{equation}
{\Phi}^{\nabla}\left(\dfrac{\partial }{\partial
{x^{i}}}\right)=X_i^H
\end{equation}
\begin{equation}
{\Phi}^{\nabla}\left({dx^j}\right)=\dfrac{\partial }{\partial
{y_{j}}}.
\end{equation}

\noindent Let $(\widehat J,\widehat g)$ be the Norden structure on  $T(M) \oplus T^{\star}(M)$ defined in previous section, the isomorphism ${\Phi}^{\nabla}$ allows us to define an almost complex structure $\tilde J$ and a neutral metric $\tilde g$ on $T^{\star}(M)$ as in the following.\\
\noindent We define $\tilde J:T(T^{\star}(M)) \rightarrow T(T^{\star}(M))$ by
\begin{equation}
{\tilde J}= ({\Phi}^{\nabla} )\circ\widehat J \circ ({\Phi}^{\nabla} )^{-1}\label{170}
\end{equation} 
\noindent and the pseudo Riemannian metric $\tilde g$ on $T^{\star}(M)$ by
\begin{equation}
{\tilde g}= (({\Phi}^{\nabla} )^{-1})^{\star}(\widehat g). \label{180}
\end{equation}

\begin{prop} $(T^{\star}(M),\tilde J,\tilde g)$ is a Norden manifold.
\end{prop}

\noindent {\bf {Proof.}} For all $X,Y\in C^{\infty}(T(T^{\star}(M)))$ we have:\\

$\tilde {g}(\tilde{J}(X),Y)=\widehat g(({\Phi}^{\nabla} )^{-1}( ({\Phi}^{\nabla} )\circ\widehat J \circ ({\Phi}^{\nabla} )^{-1})(X)),({\Phi}^{\nabla} )^{-1}(Y))$\\

$=\widehat g(\widehat J  (({\Phi}^{\nabla} )^{-1}(X)),({\Phi}^{\nabla} )^{-1}(Y))=\widehat g(({\Phi}^{\nabla} )^{-1}(X),\widehat J (({\Phi}^{\nabla} )^{-1}(Y)))$\\

$=\widehat g(({\Phi}^{\nabla} )^{-1}(X),({\Phi}^{\nabla} )^{-1}( ({\Phi}^{\nabla} )\circ \widehat J \circ  ({\Phi}^{\nabla} )^{-1})(Y)))=\tilde g(X,\tilde J (Y)). \,  \, \square $  \\

\noindent Direct computations give the following local expressions for $\tilde J$ and $\tilde g$:
\begin{equation}
\left \{
\begin{array}{l}
\vspace{0.2cm}
{\tilde J} \left(  X_i^H\right)=J^k_iX_k^H+g_{ik}{\dfrac{\partial }{\partial
y_{k}}}\\
{\tilde J} \left({\dfrac{\partial }{\partial
y_{j}}}\right)=-J^j_k{\dfrac{\partial }{\partial
y_{k}}}.
\end{array}\label{190}
\right.
\end{equation}
\begin{equation}
\left \{
\begin{array}{l}
\vspace{0.2cm}
{\tilde g} \left(  X_i^H,X^H_j\right)=g_{ij}\\
\vspace{0.2cm}
{\tilde g} \left(  X_i^H,{\dfrac{\partial }{\partial
y_{j}}}\right )={1\over 2} J^j_i\\
{\tilde g} \left({\dfrac{\partial }{\partial
y_{i}}},{\dfrac{\partial }{\partial
y_{i}}}\right)=g^{ij}.
\end{array}\label{200}
\right.
\end{equation}

\noindent Moreover, if we denote by $\tilde N$  the Nijenhuis tensor of $\tilde J$, the following hold:
\begin{equation}
\begin{array}{l}
\vspace{0.2cm}
\tilde N \left({\dfrac{\partial }{\partial y_{i}}},{\dfrac{\partial }{\partial y_{j}}}\right)={\Phi}^{\nabla}\left(N^{\nabla}(\widehat J)\left(dx^i,dx^j\right)\right) \\
\end{array}\label{210}
\end{equation}
\begin{equation}
\begin{array}{l}
\vspace{0.2cm}
\tilde N \left( X^H_{i},{\dfrac{\partial }{\partial y_{j}}}\right)=-{\left( \left( {\nabla}_{JX_{i}} J \right) X_{k}-J \left( {\nabla}_{X_{i}} J \right){ X_{k}}\right)}^j {\dfrac{\partial }{\partial y_{k}}}\\
={\Phi}^{\nabla}\left(N^{\nabla}(\widehat J)\left( {X_{i}}, dx^j\right)\right) \\
\end{array} \label{220}
\end{equation}
\begin{equation}
\begin{array}{l}
\vspace{0.2cm}
\tilde N \left(  X_i^H, X_j^H \right)={\Phi}^{\nabla}\left(N^{\nabla}(\widehat J)\left( {X_{i}},X_j\right)\right) +\\
+ y_l \left(J^k_iJ^a_jR^l_{khr} +J_h^rJ^k_iR^l_{kjr}-J^r_hJ^k_jR^l_{kir}-R^l_{ijl}\right) {{\dfrac{\partial }{\partial y_{l}}}}. \\
\end{array}\label{230}
\end{equation}

\noindent Thus we get the following.

\begin{prop} Let $(M,J,g)$ be a complex Norden manifold with flat  natural canonical connection then $(T^\star(M),\tilde J,\tilde g)$ is a complex Norden manifold.
\end{prop}

\subsection{ K\"ahler Norden structure on $T^{\star}(M)$}
\noindent A direct computation gives the following.

\begin{lem} Let $(M,J,g)$ be a Norden manifold and let $(\tilde J,\tilde g)$ be the Norden structure defined on $T^\star(M)$ as in Proposition 4.2.. Let $\nabla$ and $\tilde {\nabla}$ be the Levi Civita connection of $g$ and  $\tilde g$ respectively, then:
\begin{equation}
\begin{array}{l}
\vspace{0.2cm}
{\tilde\nabla}_{X_i^H}X_j^H=\{ {\Gamma}_{ij}^r-{\dfrac{1 }{10}}(({\nabla}_{\dfrac{\partial }{\partial x^{i}}}J){\dfrac{\partial }{\partial x^{j}}}+\\
+({\nabla}_{\dfrac{\partial }{\partial x^{j}}}J){\dfrac{\partial }{\partial x^{i}}})^p J^r_p+{\dfrac{1 }{5}}g^{rl}y_k(R^k_{ijs}J^s_l-2{R^k}_{ils}J^s_j-2R^k_{jls}J^s_i) \} X_r^H+\\
\vspace{0.2cm}
+{\dfrac{1 }{5}}\{g_{rs}(({\nabla}_{\dfrac{\partial }{\partial x^{i}}}J){\dfrac{\partial }{\partial x^{j}}}+({\nabla}_{\dfrac{\partial }{\partial x^{j}}}J){\dfrac{\partial }{\partial x^{i}}})^r+\\
+y_k(3R^k_{ijs}+J^l_s J^r_jR^k_{ilr}+J^l_sR^k_{jlr}J^r_i) \}{{\dfrac{\partial }{\partial y_{s}}}}
\end{array}\label{230}
\end{equation}
\begin{equation}
\begin{array}{l}
\vspace{0.2cm}
{\tilde\nabla}_{{{\dfrac{\partial }{\partial y_{j}}}}}X_i^H= {\dfrac{1 }{5}}\{g^{rk}(({\nabla}_{\dfrac{\partial }{\partial x^{i}}}J){\dfrac{\partial }{\partial x^{r}}}-
({\nabla}_{\dfrac{\partial }{\partial x^{r}}}J){\dfrac{\partial }{\partial x^{i}}})^j+\\
-2g^{rk}g^{jl}y_rR^r_{ijl} \} X_k^H+{\dfrac{1 }{10}}\{-J^r_s(({\nabla}_{\dfrac{\partial }{\partial x^{i}}}J){\dfrac{\partial }{\partial x^{r}}}+\\
-({\nabla}_{\dfrac{\partial }{\partial x^{r}}}J){\dfrac{\partial }{\partial x^{i}}})^s+\\
+2J^r_sg^{jl}y_kR^k_{irl} \}{{\dfrac{\partial }{\partial y_{s}}}}
\end{array}\label{250}
\end{equation}
\begin{equation}
\begin{array}{l}
\vspace{0.2cm}
{\tilde\nabla}_{X_i^H}{{{\dfrac{\partial }{\partial y_{j}}}}}={\tilde\nabla}_{{{\dfrac{\partial }{\partial y_{j}}}}}X_i^H-{\Gamma}^j_{is}{{\dfrac{\partial }{\partial y_{s}}}}
\end{array}\label{290}
\end{equation}
\begin{equation}
\begin{array}{l}
\vspace{0.2cm}
{\tilde\nabla}_{{{\dfrac{\partial }{\partial y_{i}}}}}{{{\dfrac{\partial }{\partial y_{j}}}}}=0
\end{array}\label{300}
\end{equation}
\end{lem}

\noindent In particular the following hold.
\begin{prop} Let $(M,J,g)$ be a flat  K\"ahler Norden manifold then $(T^\star(M),\tilde J,\tilde g)$ is a flat  K\"ahler Norden manifold.
\end{prop}
\noindent {\bf {Proof.}} From Lemma 4.4, under the assumption $(M,J,g)$ is K\"ahler Norden flat,  we get the following expression of $\tilde \nabla$:

\begin{equation}
\begin{array}{l}
\vspace{0.2cm}
{\tilde\nabla}_{X_i^H}X_j^H= {\Gamma}_{ij}^r X_r^H\\
\vspace{0.2cm}
\end{array}\label{310}
\end{equation}
\begin{equation}
\begin{array}{l}
\vspace{0.2cm}
{\tilde\nabla}_{{{\dfrac{\partial }{\partial y_{j}}}}}X_i^H=0
\end{array}\label{320}
\end{equation}
\begin{equation}
\begin{array}{l}
\vspace{0.2cm}
{\tilde\nabla}_{X_i^H}{{{\dfrac{\partial }{\partial y_{j}}}}}=-{\Gamma}^j_{is}{{\dfrac{\partial }{\partial y_{s}}}}
\end{array}\label{330}
\end{equation}
\begin{equation}
\begin{array}{l}
\vspace{0.2cm}
{\tilde\nabla}_{{{\dfrac{\partial }{\partial y_{i}}}}}{{{\dfrac{\partial }{\partial y_{j}}}}}=0.
\end{array}\label{340}
\end{equation}
\noindent Then we get:
\begin{equation}
\begin{array}{l}
\vspace{0.2cm}
({\tilde\nabla}_{X_i^H}{\tilde J})X_j^H=(({\nabla}_ {{{\dfrac{\partial }{\partial x^{i}}}}}J){{{\dfrac{\partial }{\partial x^{j}}}}})^r X_r^H=0\\
\vspace{0.2cm}
\end{array}\label{350}
\end{equation}
\begin{equation}
\begin{array}{l}
\vspace{0.2cm}
({\tilde\nabla}_{{{\dfrac{\partial }{\partial y_{j}}}}}{\tilde J})X_i^H=0
\end{array}\label{360}
\end{equation}
\begin{equation}
\begin{array}{l}
\vspace{0.2cm}
({\tilde\nabla}_{X_i^H}{\tilde J}){{{\dfrac{\partial }{\partial y_{j}}}}}=-(({\nabla}_ {{{\dfrac{\partial }{\partial x^{i}}}}}J){{{\dfrac{\partial }{\partial x^{r}}}}})^j{{{\dfrac{\partial }{\partial y_{j}}}}} =0\\
\end{array}\label{370}
\end{equation}
\begin{equation}
\begin{array}{l}
\vspace{0.2cm}
({\tilde\nabla}_{{{\dfrac{\partial }{\partial y_{i}}}}}{\tilde J}){{{\dfrac{\partial }{\partial y_{j}}}}}=0.
\end{array}\label{380}
\end{equation}
\noindent Hence $\tilde \nabla \tilde J=0$.\\

\noindent Moreover let $\tilde R$ be the Riemann curvature tensor of $\tilde g$, we have:
\begin{equation}
\begin{array}{l}
\vspace{0.2cm}
{\tilde R}({X_i^H}, X_j^H)X_k^H=R^l_{ijk} X^H_l=0\\
\vspace{0.2cm}
\end{array}\label{390}
\end{equation}
\begin{equation}
\begin{array}{l}
\vspace{0.2cm}
{\tilde R}({X_i^H}, X_j^H){{{\dfrac{\partial }{\partial y_{k}}}}}=R^k_{jir}{{{\dfrac{\partial }{\partial y_{r}}}}}=0
\end{array}\label{400}
\end{equation}
\begin{equation}
\begin{array}{l}
\vspace{0.2cm}
{\tilde R}({X_i^H},{{{\dfrac{\partial }{\partial y_{k}}}}}) X_j^H=0\\
\end{array}\label{410}
\end{equation}
\begin{equation}
\begin{array}{l}
\vspace{0.2cm}
{\tilde R}({X_i^H},{{{\dfrac{\partial }{\partial y_{j}}}}}) {{{\dfrac{\partial }{\partial y_{k}}}}}=0\\
\end{array}\label{420}
\end{equation}
\begin{equation}
\begin{array}{l}
\vspace{0.2cm}
{\tilde R}({{{\dfrac{\partial }{\partial y_{i}}}}},{{{\dfrac{\partial }{\partial y_{j}}}}}) {{{\dfrac{\partial }{\partial y_{k}}}}}=0.\\
\end{array}\label{430}
\end{equation}
\noindent Thus $\tilde R=0$ and the proof is complete. $ \square $

\vspace{0.1cm}
\noindent {\bf author address}: Dipartimento di Matematica e Informatica "U. Dini" Viale Morgagni 67/a 50134 Firenze, Italy, \, {\bf email}: antonella.nannicini@unifi.it

\end{document}